\documentclass[9pt,shortpaper,twoside,web]{ieeecolor}

\usepackage{generic}
\usepackage{cite}
\usepackage{amsmath,amssymb,amsfonts}
\usepackage{algorithmic}
\usepackage{graphicx}
\usepackage{nicematrix}
\usepackage{tikz,mathrsfs}
\usepackage{subfigure}
\usetikzlibrary{fit}
\usepackage{mathtools}

\usepackage{textcomp}
\def\BibTeX{{\rm B\kern-.05em{\sc i\kern-.025em b}\kern-.08em
		T\kern-.1667em\lower.7ex\hbox{E}\kern-.125emX}}
\newtheorem{theorem}{Theorem}

\newtheorem{lemma}{Lemma}

\newtheorem{remark}{Remark}

\newtheorem{definition}{Definition}

\DeclareMathOperator{\esssup}{ess~sup}
\DeclareMathOperator{\blkdiag}{blk-diag}

\DeclareMathOperator{\rank}{rank}
\DeclareMathOperator{\loc}{loc}
\DeclareMathOperator{\RE}{Re}

\newcommand{\R}{{\mathbb{R}}}

\newcommand\myae{\stackrel{\mathclap{\normalfont\mbox{a.e.}}}{=}}
\newcommand{\setdef}[2]{\left\{\, #1 \left|\, \vphantom{#1} #2\right.\right\}}

\newenvironment{smallbmatrix}
{\left[\begin{smallmatrix}}
{\end{smallmatrix}\right]}

\begin{document}
	\title{Partial detectability and generalized functional observer design for linear descriptor systems}
	\author{Juhi Jaiswal, Thomas Berger, and Nutan Kumar Tomar$^{*}$
		\thanks{Manuscript received January $24$, $2023$.}
		\thanks{This work is supported by the Science and Engineering Research Board, New Delhi through the Grant MTR/2019/000494.}
		\thanks{Juhi Jaiswal and Nutan Kumar Tomar are with the Department of Mathematics, Indian Institute of Technology Patna, India (e-mail: juhi$\_$1821ma03@iitp.ac.in, nktomar@iitp.ac.in).}%
		\thanks{Thomas Berger is with Universit\"at Paderborn, Institut f\"ur Mathematik, Warburger Str.~100, 33098~Paderborn, Germany (e-mail: thomas.berger@math.upb.de).}%
		\thanks{$^{*}$ Corresponding author.}}

	\maketitle
	
	\begin{abstract}
		This paper studies linear time-invariant descriptor systems which are not necessarily regular. We introduce the notion of partial detectability and characterize this concept by means of a simple rank criterion involving the system coefficient matrices. Three particular cases of this characterization are discussed in detail. Furthermore, we show that partial detectability is necessary for the existence of a generalized functional observer, but not sufficient. We identify a condition which together with partial detectability gives sufficiency.
	\end{abstract}
	
	\begin{IEEEkeywords}
		Descriptor systems, Partial detectability, Wong sequences, Observer design, Generalized functional observer
	\end{IEEEkeywords}
	
\section{Introduction}\label{intro}
We consider linear time-invariant (LTI) descriptor systems of the form
\begin{subequations}\label{dls}
	\begin{IEEEeqnarray}{rCl}
		E\dot{x}(t) &=& Ax(t) + Bu(t),  \label{dlsa} \\
		y(t) &=& Cx(t) + Du(t) , \label{dlsb} \\
		z(t) &=& Kx(t), \label{dlsc}
	\end{IEEEeqnarray}
\end{subequations}	where $E,~A \in \mathbb{R}^{m \times n},~B \in \mathbb{R}^{m \times k},~C \in \mathbb{R}^{p \times n},D \in \mathbb{R}^{p \times k}$, and $K \in \mathbb{R}^{r \times n}$ are known constant matrices. Systems of type \eqref{dls} are also called singular systems or systems described by differential-algebraic equations (DAEs). The first order matrix polynomial $(\lambda E-A)$, in the indeterminate $\lambda$, is called matrix pencil for \eqref{dls}. System~\eqref{dls} is called regular if $m=n$ and $\det (\lambda E-A)$ is not the zero polynomial in~$\lambda$. In the present paper, we do not assume that the system is regular; in fact, no assumptions on the coefficient matrices in~\eqref{dls} are made, and the system may be under- and/or over-determined. Such systems occur naturally when dynamical systems are subject to algebraic constraints; for further motivation, see the books~\cite{campbell1982singular,dai1989singular,kunkel2006differential,riaza2008differential} and the references therein.
	
	We call $x: \mathbb{R} \rightarrow \mathbb{R}^n$ the semistate of the system~\eqref{dls}, because unlike state space systems, $x(t)$ does not satisfy the semigroup property and cannot be arbitrarily initialized~\cite{newcomb1981semistate}. However, $x(t)$ contains the full information about all intrinsic properties of the system at time $t$. The functions $u: \mathbb{R} \rightarrow \mathbb{R}^k$ and $y: \mathbb{R} \rightarrow \mathbb{R}^p$ are called the input and the output of system~\eqref{dls}, respectively, and they are obtained by measurements, e.g., via sensors. The functional vector $z(t) \in \mathbb{R}^r$ contains those variables which cannot be measured, and observers are required to estimate them, cf.\ Definition~\ref{def:observer}. If $K$ is not the identity matrix, such an observer is called a functional (or partial state) observer; otherwise, we call it a full state observer. The tuple $(x,u,y,z):\mathbb{R} \rightarrow \mathbb{R}^n \times \mathbb{R}^k \times \mathbb{R}^p \times \mathbb{R}^r$ is said to be a solution of~\eqref{dls}, if it belongs to the set
	\begin{IEEEeqnarray*}{rCl}
		\mathscr{B} &:=& \{(x,u,y,z) \!\in\! \mathscr{L}^1_{\loc}(\mathbb{R}; \mathbb{R}^{n+k+p+r}) \mid Ex \in \mathcal{AC}_{\loc}(\mathbb{R} ; \mathbb{R}^m)  \\
		&&\hfill\raggedright   \text{ and } (x,u,y,z)   \text{ satisfies } \eqref{dls} \text{ for almost all } t \in \mathbb{R} \},
	\end{IEEEeqnarray*}
	where $\mathscr{L}^1_{\loc}$ is the set of measurable and locally Lebesgue integrable functions and $\mathcal{AC}_{\loc}$ represents the set of locally absolutely continuous functions. Descriptor systems based on the \textit{behavior} $\mathscr{B}$ have been studied in detail e.g.\ in~\cite{berger2014}. Exploiting the behavior $\mathscr{B}$, various observability and detectability concepts for descriptor systems~\eqref{dls} are studied in~\cite{berger2017observability}. Existence conditions for full state and functional observers of descriptor systems have been investigated in~\cite{berger2019ode,jaiswal2021necessary,berger2017observers,berger2019disturbance,jaiswal2021existence,jaiswal2021functional,jaiswal2022icc}, see also the references therein. Throughout the article, we assume that the behavior $\mathscr{B}$ is nonempty, which amounts to the existence of an admissible pair for~\eqref{dls}, consisting of an admissible initial condition and input function, see also~\cite{jaiswal2021necessary}.
	
	The major contribution of the current study is twofold. First, we prove that the system~\eqref{dls} is partially detectable if, and only if, a simple rank condition involving the system coefficient matrices holds. This result first requires to establish a precise definition of partial detectability for~\eqref{dls} based on the behavior~$\mathscr{B}$. Roughly speaking, detectability means that the inputs and outputs determine the state variables asymptotically. From this point of view, the partial detectability of~\eqref{dls} is related to the asymptotic determination of~$z(t)$ from the knowledge of~$u(t)$ and~$y(t)$. In a particular case of the characterization of partial detectability, we deduce that the existing algebraic characterization of partial detectability of state space systems in~\cite[Thm.~1]{darouach2022functional} may give an erroneous result.  The second major contribution is to relate the concept of partial detectability to the existence of generalized functional observers. We show that partial detectability is necessary for their existence but not sufficient. Furthermore, we derive an additional condition, which, together with partial detectability, allows for the construction of a generalized functional observer. For this, we provide a step-by-step algorithm. Our approach is purely algebraic and based on simple matrix theory.
	
	The paper is organized as follows. Section \ref{prelim} collects some preliminary results used in the sequel of the article. In Section~\ref{mainresult}, the concept of partial detectability of system~\eqref{dls} is introduced along with algebraic test conditions and their equivalence is shown.  Section~\ref{particularcases} discusses some particular cases of the proposed results and emphasizes a significant modification to the existing theory of partial detectability for state space systems.  In Section~\ref{obsvdesign}, we show that partial detectability is necessary for the existence of a generalized functional observer for system~\eqref{dls} and identify a condition, together with which it is also sufficient. Section~\ref{numerical} considers a numerical example to illustrate the design algorithm. Finally, Section~\ref{conc} concludes the article.
	
	We use the following notations: $0$ and $I$ stand for appropriate dimensional zero and identity matrices, respectively. For more clarity, the identity matrix of size $n \times n$ is sometimes denoted by $I_n$. The set of complex numbers is denoted by $\mathbb{C}$, $\bar{\mathbb{C}}^+ := \{\lambda \in \mathbb{C}~| ~ \RE(\lambda)\geq 0\}$, and $\mathbb{C}^- := \{\lambda \in \mathbb{C}~| ~ \RE(\lambda)< 0\}$. The symbols $A^\top$ and $\ker A$ denote the transpose and null space of a matrix $A \in \mathbb{R}^{m \times n}$, respectively. In a block partitioned matrix, all missing blocks are zero matrices of appropriate dimensions. The set $\sigma(M)$ denotes the spectrum of a matrix $M \in \mathbb{R}^{n \times n}$.
	The notation $f \myae g$ means $f,~g \in \mathscr{L}^1_{\loc}(\mathbb{R};\mathbb{R}^n)$ are equal almost everywhere, \emph{i.e.}, $f(t) = g(t)$ for almost all $t \in \mathbb{R}$. The set $AM := \setdef{Ax}{x \in M}$ is the image of a subspace $M \subseteq \mathbb{R}^n$ under $A \in \mathbb{R}^{m \times n}$ and $A^{-1}M := \setdef{x \in \mathbb{R}^n}{Ax \in M}$ represents the pre-image of $M \subseteq \mathbb{R}^m$ under $A \in \mathbb{R}^{m \times n}$. A block matrix having diagonal elements $A_1, \ldots, A_k$ is represented by $\blkdiag(A_1,\ldots,A_k)$.
	
	\section{Preliminaries}\label{prelim}
	In this section, we first recall some preliminary results from  matrix theory and the theory of descriptor systems. These are fundamental to the development of the main results in this paper.	
	
\begin{lemma}\cite[Quasi-Kronecker Form (QKF)]{berger2013addition}\label{pre:lm1}
		For every matrix pencil $(\lambda E - A)$ with $E,A\in\R^{m\times n}$ there exist nonsingular matrices $P \in \mathbb{R}^{m \times m}$ and $Q \in \mathbb{R}^{n \times n}$ such that
		\begin{equation}\label{eq:qkf}
		P(\lambda E - A)Q =  \begin{bmatrix}
		\lambda E_{\epsilon} - A_{\epsilon} & & & \\ & \lambda I_f - J_f & & \\ & & \lambda J_{\sigma} - I_{\sigma} & \\ & & & \lambda E_{\eta} - A_{\eta}
		\end{bmatrix}
		\end{equation}
		where
\begin{enumerate}
	\item $E_{\epsilon} ,~ A_{\epsilon} \in \mathbb{R}^{m_{\epsilon} \times n_{\epsilon}}$, $m_{\epsilon} < n_{\epsilon}$, and $\rank E_{\epsilon} = \rank (\lambda E_{\epsilon} - A_{\epsilon}) = m_{\epsilon}$ for all $\lambda \in \mathbb{C}$.
	\item $J_f \in \mathbb{R}^{n_f \times n_f}$.
	\item $J_{\sigma} \in \mathbb{R}^{n_{\sigma} \times n_{\sigma}}$ is nilpotent.
	\item $E_{\eta} ,~ A_{\eta} \in \mathbb{R}^{m_{\eta} \times n_{\eta}}$, $m_{\eta} > n_{\eta}$, and $\rank E_{\eta} = \rank(\lambda E_{\eta} - A_{\eta}) = n_{\eta} $ for all $\lambda \in \mathbb{C}$.
\end{enumerate}
\end{lemma}
	
	
	\begin{lemma}\label{pre:lm5}
		For any matrices $X$ and $Y$of compatible dimensions, $\rank \begin{bmatrix}
		X \\ Y
		\end{bmatrix} = \rank X$ if, and only if, $\ker{X} \subseteq \ker{Y}$.
	\end{lemma}
	
	\begin{lemma}\cite{jaiswal2021necessary}\label{pre:lm2}
		Let $X \in \mathbb{R}^{m_1 \times r_1},~S \in \mathbb{R}^{m_1 \times r_2}$, and $Y \in \mathbb{R}^{m_2 \times r_2}$. If $X$ has full row rank and/or $Y$ has full column rank, then
		$$\rank \begin{bmatrix}
		X & S \\ 0 & Y
		\end{bmatrix} = \rank{X} + \rank{Y}.$$
	\end{lemma}
	
	We conclude this section by recalling the concept of a complex Wong sequence corresponding to \eqref{dlsa} from~\cite{berger2019disturbance}.
	\begin{definition}\label{def:wong:detectability}
		For a given system \eqref{dlsa} with $E,A\in\R^{m\times n}$ and $\lambda\in\mathbb{C}$ the Wong sequence $\{\mathcal{W}_{[E,A],\lambda}^i\}_{i = 0}^{\infty}$ is a sequence of complex subspaces, defined by
		\begin{eqnarray*}
			&& \mathcal{W}_{[E,A],\lambda}^0 := \{0\}, \\
			&&  \mathcal{W}_{[E,A],\lambda}^{i+1} := (A- \lambda E)^{-1}(E\mathcal{W}_{[E,A],\lambda}^i) \subseteq \mathbb{C}^n, \\
			&& \mathcal{W}_{[E,A],\lambda}^* := \bigcup_{i \in \mathbb{N}} \mathcal{W}_{[E,A],\lambda}^i .
		\end{eqnarray*}
	\end{definition}
	
	\section{Partial detectability}\label{mainresult}
	
The main aim of this section is to derive a simple rank condition for partial detectability of~\eqref{dls} in terms of the system coefficient matrices. First, we define the concept of partial detectability of~\eqref{dls} in terms of the behavior~$\mathscr{B}$, which is a natural extension of the detectability of \eqref{dlsa}--\eqref{dlsb}. Throughout the paper, the notation $``x(t) \rightarrow 0 \text{ as } t \rightarrow \infty"$  means $``\lim\limits_{t \rightarrow \infty} \esssup\limits_{ [t, \infty)} ||x(t)|| = 0"$.
	
	\begin{definition}\label{def:kdetectability}
		The descriptor system~\eqref{dls} is said to be partially detectable, if for all $(x_1,u,y,z_1), (x_2,u,y,z_2)\in \mathscr{B}$ we have that  $z_1(t) - z_2(t) \rightarrow 0$ as $t \rightarrow \infty$.
	\end{definition}

By linearity of the behavior~$\mathscr{B}$ it is clear that~\eqref{dls} is partially detectable if, and only if, for all $(x,0,0,z)\in \mathscr{B}$ we have that  $z(t) \rightarrow 0$ as $t \rightarrow \infty$. Therefore, partial detectability is independent of the matrices~$B$ and~$D$.
	
\begin{remark}\label{rem2}
In \eqref{dls}, if $K=I_n$, the above definition reduces to the detectability of \eqref{dlsa}--\eqref{dlsb}, see \cite{hou1999observer,gupta2016detectability,berger2017observability,jaiswal2021necessary}. Note that detectability is called ``behavioral detectability'' in \cite{berger2017observability}.
	\end{remark}
	
The aim of the remainder of this section is to derive a rank criterion for partial detectability of~\eqref{dls}. Additionally to the new rank criterion, we include a characterization of partial detectability  in terms of the Wong sequences, which was already implicitly contained in \cite{berger2019disturbance}.

Let $l \in \mathbb{N}$, $\lambda\in\mathbb{C}$ and introduce the following notations: 	
	\begin{IEEEeqnarray*}{rCl}
		&	\mathcal{E} = \begin{bmatrix}
			E \\ 0 \end{bmatrix},\quad \mathcal{A} = \begin{bmatrix} A \\ C \end{bmatrix}, \\
		& \mathcal{G}_{l,[\mathcal{E},\mathcal{A}],\lambda} =  \NiceMatrixOptions
		{nullify-dots,code-for-last-col = \color{black},code-for-last-col=\color{black}}
		\begin{bNiceMatrix}[first-row,last-col]
			& \Ldots[line-style={solid,<->},shorten=0pt]^{l \text{ block columns}} \\
			\lambda \mathcal{E} - \mathcal{A} & & & & &  \\
			\mathcal{E} & \lambda \mathcal{E} - \mathcal{A} & & & &  \\
			 & \ddots & \ddots &   &  &   \\
			& &  \mathcal{E} & \lambda \mathcal{E} - \mathcal{A} & &  \\
			& &  &\mathcal{E} & \lambda \mathcal{E} - \mathcal{A} &
		\end{bNiceMatrix},&
	\end{IEEEeqnarray*}
	and
	\begin{IEEEeqnarray*}{rCl}
		& \mathcal{G}_{l,[\mathcal{E},\mathcal{A},K],\lambda} =  \NiceMatrixOptions
		{nullify-dots,code-for-last-col = \color{black},code-for-last-col=\color{black}}
		\begin{bNiceMatrix}[first-row,last-col]
			& \Ldots[line-style={solid,<->},shorten=0pt]^{l \text{ block columns}} \\
			\lambda \mathcal{E} - \mathcal{A} & & & & &   \\
			\mathcal{E} & \lambda \mathcal{E} - \mathcal{A} & & & &  \\
			& \ddots & \ddots &   &  &  \\
			& &  \mathcal{E} & \lambda \mathcal{E} - \mathcal{A} & &  \\
			& &  &\mathcal{E} & \lambda \mathcal{E} - \mathcal{A} & \\
			& &  & & K &
		\end{bNiceMatrix}. &
	\end{IEEEeqnarray*}

We are now ready to state the first main result of this paper.
	
\begin{theorem}\label{thm:kdetectability}
For a given system \eqref{dls} the following statements are equivalent:
\begin{enumerate}
	\item[(a)]  The system \eqref{dls} is partially detectable. \label{a}
	\item[(b)] The following condition holds: \label{b}
	\begin{equation}\label{kdetectability}
	\forall\,\lambda \in \bar{\mathbb{C}}^+:\ \rank \mathcal{G}_{n,[\mathcal{E},\mathcal{A}],\lambda} = \rank \mathcal{G}_{n,[\mathcal{E},\mathcal{A},K],\lambda}.
			\end{equation}
	\item[(c)] $\forall\,\lambda \in \bar{\mathbb{C}}^+:\ \mathcal{W}_{[\mathcal{E},\mathcal{A}],\lambda}^{*} \subseteq \ker{K}$. \label{c}
		\end{enumerate}
\end{theorem}
	
\begin{proof}
	(a) $\Leftrightarrow$ (c): This follows from \cite[Lem.~A.4]{berger2019disturbance}.
		
	(b) $\Rightarrow$ (c): Fix $\lambda \in \bar{\mathbb{C}}^+$. In view of Lemma~\ref{pre:lm5}, condition (b) is equivalent to
		\begin{IEEEeqnarray}{rCl}
			& \ker \mathcal{G}_{n,[\mathcal{E},\mathcal{A}],\lambda} \subseteq \ker
			\NiceMatrixOptions
			{nullify-dots,code-for-last-col = \color{black},code-for-last-col=\color{black}}
			\begin{bNiceMatrix}[first-row]
				& \Ldots[line-style={solid,<->},shorten=0pt]^{n \text{ block columns}} \\
				0 & 0 & \hdots & 0 & K  \end{bNiceMatrix}.& \label{eq:proof1}
		\end{IEEEeqnarray}
Now let $x\in \mathcal{W}_{[\mathcal{E},\mathcal{A}],\lambda}^{*}$ and observe that, since the Wong sequence $\{\mathcal{W}_{[\mathcal{E},\mathcal{A}],\lambda}^{i}\}_{i=0}^\infty$ terminates after finitely many steps, and in each step before termination the dimension of the associated space increases by at least one, it is clear that  $\mathcal{W}_{[\mathcal{E},\mathcal{A}],\lambda}^{n} = \mathcal{W}_{[\mathcal{E},\mathcal{A}],\lambda}^{*}$. Thus $x\in \mathcal{W}_{[\mathcal{E},\mathcal{A}],\lambda}^n$ and there exist $x_1,\ldots,x_n \in \mathbb{C}^n$ such that $x=x_n$ and
\begin{IEEEeqnarray}{rCl}
		& (\lambda  \mathcal{E} -  \mathcal{A}) x_1 = 0 ,~
			\mathcal{E}x_1 + (\lambda  \mathcal{E} -  \mathcal{A}) x_2 = 0,~ \ldots,& \notag \\ &
			\qquad \qquad \mathcal{E}x_{n-1} + (\lambda  \mathcal{E} -  \mathcal{A}) x_n  = 0.
			 & \label{eq:proof2}
\end{IEEEeqnarray}
Therefore, $(x_1^\top, \ldots,x_n^\top)^\top \in \ker \mathcal{G}_{n,[\mathcal{E},\mathcal{A}],\lambda}$ and by~\eqref{eq:proof1} this gives $x = x_n \in \ker K$.
		
(c) $\Rightarrow$ (b): Fix $\lambda \in \bar{\mathbb{C}}^+$ and let $x\in \ker \mathcal{G}_{n,[\mathcal{E},\mathcal{A}],\lambda}$. We show that $x\in \ker [0,\ldots, 0,K]$, which proves~\eqref{eq:proof1} and hence also~(b). Write $x = (x_1^\top, \ldots,x_n^\top)^\top$ with $x_1,\ldots,x_n \in \mathbb{C}^n$, then \eqref{eq:proof2} holds.

Therefore, $x_i \in \mathcal{W}_{[\mathcal{E},\mathcal{A}],\lambda}^i$ for $i=1,\ldots,n$ and since $\mathcal{W}_{[\mathcal{E},\mathcal{A}],\lambda}^{*} = \mathcal{W}_{[\mathcal{E},\mathcal{A}],\lambda}^n$ it follows from~(c) that $x_n \in \ker K$. This completes the proof.
	\end{proof}
	
In the following lemma, we derive a characterization for partial detectability in terms of the QKF of the matrix pencil $(\lambda \mathcal{E} - \mathcal{A})$. Later, this will be used in Section~\ref{obsvdesign} to design a generalized functional observer.

		\begin{lemma}\label{lm:obsv}
		Consider a system~\eqref{dls} and let the matrix pencil $(\lambda \mathcal{E} - \mathcal{A})$ have QKF~\eqref{eq:qkf} such that $K$ is partitioned accordingly, i.e.,
        \[
            K Q = \begin{bmatrix}	K_{\epsilon} & K_f & K_{\sigma} & K_{\eta}  \end{bmatrix}.
        \]
        Then~\eqref{dls} is partially detectable if, and only if,
		\begin{enumerate}
		    \item[(a)] $K_{\epsilon} = 0$ and
		    \item[(b)] $\forall\,\lambda \in \bar{\mathbb{C}}^+:\ \ker (\lambda I_f - J_f)^n \subseteq \ker K_f $.
		\end{enumerate}
		Moreover, if $J_f = \begin{bmatrix}
		J_{f_1} \\ & J_{f_2}
		\end{bmatrix}$ and $K_f = \begin{bmatrix}
		K_{f_1} & K_{f_2}
		\end{bmatrix}$ such that $\sigma(J_{f_1}) \subseteq \bar{\mathbb{C}}^+$ and $\sigma(J_{f_2}) \subseteq \mathbb{C}^-$, then (b) is equivalent to $K_{f_1} = 0$.
	\end{lemma}
	
	\begin{proof}
        Without loss of generality, assume that $(\lambda \mathcal{E} - \mathcal{A})$ is in QKF~\eqref{eq:qkf} and
        \[
            K = \begin{bmatrix}	K_{\epsilon} & K_f & K_{\sigma} & K_{\eta}  \end{bmatrix}.
        \]
        Then, as shown in~\cite[Eq.~(A.4)]{berger2019disturbance}, we have that
        \[
            \mathcal{W}_{[\mathcal{E},\mathcal{A}],\lambda}^{*}  = \mathbb{C}^{n_{\epsilon}} \times \ker (\lambda I_f - J_f)^n \times \{0\}^{n_\sigma} \times \{0\}^{n_\eta}.
        \]
		By condition~(c) in Theorem \ref{thm:kdetectability} it follows that partial detectability is equivalent to
		\begin{equation*}
		\mathbb{C}^{n_{\epsilon}} \subseteq \ker K_{\epsilon} \quad\text{and}\quad  \ker (\lambda I_f - J_f)^n  \subseteq \ker  K_f
		\end{equation*}
        for all $\lambda \in \bar{\mathbb{C}}^+$, which is equivalent to~(a) and~(b).
		
		If $J_f = \begin{bmatrix}
		J_{f_1} \\ & J_{f_2}
		\end{bmatrix}$ and $K_f = \begin{bmatrix}
		K_{f_1} & K_{f_2}
		\end{bmatrix}$, where $\lambda(J_{f_1}) \subseteq \bar{\mathbb{C}}^+$ and $\lambda(J_{f_2}) \subseteq \mathbb{C}^-$, then~(b) is equivalent to
		\begin{IEEEeqnarray*}{rCl}
			\forall\,\lambda \in \bar{\mathbb{C}}^+:\ \ker \begin{bmatrix}
		(\lambda I \!-\! J_{f_1})^n \\ & (\lambda I \!-\! J_{f_2})^n
		\end{bmatrix} \subseteq \ker  \begin{bmatrix} K_{f_1} & K_{f_2}	\end{bmatrix}.
		\end{IEEEeqnarray*}
    Since $\sigma(J_{f_2}) \subseteq \mathbb{C}^-$ we have $\ker (\lambda I - J_{f_2})^n = \{0\}$, so the above condition is equivalent to
		\begin{IEEEeqnarray*}{rCl}
		\forall\,\lambda \in \bar{\mathbb{C}}^+:\ \ker(\lambda I - J_{f_1})^n  &\subseteq& \ker K_{f_1}.
	\end{IEEEeqnarray*}
Since $\sigma(J_{f_1}) \subseteq \bar{\mathbb{C}}^+$, this condition is in turn equivalent to
\begin{IEEEeqnarray*}{rCl}\label{eq:f1}
  \mathbb{C}^{n_{f_1}} = \bigcup_{\lambda \in \bar{\mathbb{C}}^+} \ker(\lambda I - J_{f_1})^n  &\subseteq& \ker K_{f_1},
\end{IEEEeqnarray*}
	where $n_{f_1}$ is the dimension of the square matrix $J_{f_1}$, thus $K_{f_1} = 0$.
	\end{proof}

The following remark is warranted on Theorem~\ref{thm:kdetectability} and Lemma~\ref{lm:obsv}.

\begin{remark}\label{rem3}
It is apparent from the proof of Theorem \ref{thm:kdetectability} that, if
	$\mathcal{W}_{[\mathcal{E},\mathcal{A}],\lambda}^{*} = \mathcal{W}_{[\mathcal{E},\mathcal{A}],\lambda}^{s}$ for some $s \in \mathbb{N}$, then the number $n$ in statement (b) of Theorem \ref{thm:kdetectability} can be replaced by $s$. Here, we use $n$ because $s$ is not known in advance and our main aim is to provide a condition directly in terms of the known data, \emph{i.e.}, the system coefficient matrices. Moreover, to verify the partial detectability of~\eqref{dls}, it is sufficient to check condition \eqref{kdetectability} in Theorem \ref{thm:kdetectability} only for those finite eigenvalues of the matrix pencil $( \lambda \mathcal{E} - \mathcal{A})$ which belong to $\bar{\mathbb{C}}^+$. For the computation of finite eigenvalues, it is recommended to use the QKF. 
\end{remark}

	\section{Particular cases}\label{particularcases}
	In this section, we discuss three particular cases of Theorem~\ref{thm:kdetectability}.
	
	\subsection{Detectability of descriptor system \eqref{dlsa}--\eqref{dlsb}}
	This case corresponds to $K = I_n$ in \eqref{dls}. By substituting $K = I_n$ in statement (c) in Theorem \ref{thm:kdetectability}, it reduces to $\mathcal{W}_{[\mathcal{E},\mathcal{A}],\lambda}^{*} = \{0\}$, which means that, for each $i \in \mathbb{N}$ and $\lambda \in \bar{\mathbb{C}}^+$ we have $\mathcal{W}_{[\mathcal{E},\mathcal{A}],\lambda}^{i} = \{0\}$.
	For $i=1$ we obtain that, for all $\lambda \in \bar{\mathbb{C}}^+$, $\ker(\lambda\mathcal{E}-\mathcal{A}) = \{0\}$, \emph{i.e.},  \begin{equation}\label{detectability}
	\forall\,\lambda \in \bar{\mathbb{C}}^+:\ \rank \begin{bmatrix}
	\lambda E -A \\ C
	\end{bmatrix} = n,
	\end{equation} which is the standard characterization of detectability of \eqref{dlsa}--\eqref{dlsb}.
	
	\subsection{Partial detectability of state space systems}
	If $E = I_n$ in \eqref{dls}, then Definition~\ref{def:kdetectability} provides the notion of partial detectability of standard state space systems. Moreover, the algebraic criterion~\eqref{kdetectability} reduces to the condition
	\begin{equation}\label{ss:kdetectability}
	\forall\,\lambda \in \bar{\mathbb{C}}^+:\ \rank \begin{bmatrix}
	(\lambda I-A)^n \\ C(\lambda I-A)^{n-1} \\ \vdots \\ C(\lambda I-A) \\ C \\ K \end{bmatrix} =
	\rank \begin{bmatrix}
	(\lambda I-A)^n \\C(\lambda I-A)^{n-1} \\ \vdots \\ C(\lambda I-A) \\ C \end{bmatrix}.
	\end{equation}
	To prove \eqref{ss:kdetectability}, first consider $\mathcal{G}_{2,[\mathcal{E},\mathcal{A}],\lambda}$ and substitute $\mathcal{E} = \begin{bmatrix} I \\ 0 \end{bmatrix}$ and $\mathcal{A} = \begin{bmatrix} A \\ C
	\end{bmatrix}$, thus obtaining
	\begin{equation}\label{n1}
	\mathcal{G}_{2,[\mathcal{E},\mathcal{A}],\lambda} = \begin{bmatrix}
	\lambda I - A & 0 \\ -C & 0 \\ I & \lambda I - A \\ 0 & -C \end{bmatrix}.
	\end{equation}
	Then, by multiplying the matrix on the right hand side of~\eqref{n1} with $U = \begin{bmatrix}
	I_n  && -(\lambda I-A) \\ & I & C \\ && I_n \\ &&& I
	\end{bmatrix}$ from the left, and then applying Lemma~\ref{pre:lm2} with $X=I_n$, we obtain
	\begin{equation*}
	\rank \mathcal{G}_{2,[\mathcal{E},\mathcal{A}],\lambda} =
	n + \rank \begin{bmatrix}
	(\lambda I - A)^2 \\ C(\lambda I-A) \\ C
	\end{bmatrix}.
	\end{equation*}
	By a similar calculation, it follows that
	\begin{equation*}
	\rank \mathcal{G}_{2,[\mathcal{E},\mathcal{A},K],\lambda} = n + \rank \begin{bmatrix}
	(\lambda I - A)^2 \\ C(\lambda I-A) \\ C \\ K
	\end{bmatrix}.
	\end{equation*}
	With this argument applied to  $\mathcal{G}_{n,[\mathcal{E},\mathcal{A}],\lambda}$ and $\mathcal{G}_{n,[\mathcal{E},\mathcal{A},K],\lambda}$ for $n>2$, it can be shown that~\eqref{kdetectability} reduces to~\eqref{ss:kdetectability}.
	
	In the articles~\cite{fernando2010functional,mohajerpoor2016new,darouach2022functional}, it has been reported (see e.g.~\cite[Thm.~1]{darouach2022functional}) that a state space system is partially detectable (note that the notion is called ``functional detectability'' in these works) if, and only if,
	\begin{equation}\label{ss:kdetectability1}
	\forall\,\lambda \in \bar{\mathbb{C}}^+:\ \rank \begin{bmatrix} \lambda I-A \\ C \\ K \end{bmatrix} =
	\rank \begin{bmatrix} \lambda I-A \\ C \end{bmatrix}.
	\end{equation}
However, this condition is obviously not equivalent to~\eqref{ss:kdetectability} and is incorrect in general. In the aforementioned works, it has been implicitly used that all eigenvalues of the matrix $A$ with nonnegative real part are semisimple (i.e., their algebraic and geometric multiplicities coincide), but this assumption was not stated explicitly. In fact, if all eigenvalues of $A$ having nonnegative real part are semisimple, then it is not hard to show that~\eqref{ss:kdetectability} reduces to~\eqref{ss:kdetectability1}. As an explicit counterexample for condition~\eqref{ss:kdetectability1}, consider the state space system $\dot{x}(t) = Ax(t)$, $y(t) = Cx(t)$, $z(t) = Kx(t)$, with
	\begin{IEEEeqnarray*}{rCl}
		A= \begin{bmatrix}  1 & 1 \\  0 & 1 \end{bmatrix},~
		C = \begin{bmatrix} 0 & 0  \end{bmatrix},
		\text{ and } K = \begin{bmatrix} 0 &  1 \end{bmatrix}.
	\end{IEEEeqnarray*}
	It is easy to verify that~\eqref{ss:kdetectability1} is satisfied for $\lambda = 1$ and hence for all $\lambda \in \bar{\mathbb{C}}^+$. On the other hand,  for $\lambda = 1$ condition~\eqref{ss:kdetectability} is not satisfied. And indeed, the system is not partially detectable: For initial data $x(0) = \begin{bmatrix}
	x_1^0 & x_2^0  \end{bmatrix}^\top$ the solution is given by
	\begin{IEEEeqnarray*}{rCl}
		x(t) &=& \exp(At)x(0) = \begin{bmatrix}  \exp(t) & t\exp(t) \\  0 & \exp(t) \end{bmatrix}\begin{bmatrix}x_1^0 \\ x_2^0 \end{bmatrix}, \\
		y(t) &=& 0, \quad
		z(t) = x_2(t) = x_2^0 \exp(t).
	\end{IEEEeqnarray*}
	Thus, for $x_2^0 \neq 0$ we have $y = 0$, but $z(t) \nrightarrow 0$ as $t \rightarrow \infty$.
	
	\subsection{Detectability of state space systems}
	Let $E=I_n$ and $K=I_n$ in \eqref{dls}. Then it is clear from~\eqref{detectability} that~\eqref{kdetectability} reduces to
	\begin{equation*}\label{ss:detectability1}
	\forall\,\lambda \in \bar{\mathbb{C}}^+:\ \rank \begin{bmatrix} \lambda I-A \\ C \end{bmatrix} =
	n,
	\end{equation*}
	which is the classical characterization of detectability for standard state space systems.

\section{Observer design}\label{obsvdesign}
In this section, we propose an observer of the following form to estimate the functional vector $z(t)$ in \eqref{dls}:
	\begin{subequations}{\label{obs}}
		\begin{IEEEeqnarray}{rCl}
			\dot{w}(t) &=& Nw(t) + H \begin{bmatrix}
				u(t) \\ y(t)
			\end{bmatrix},  \label{obsa} \\
			\hat{z}(t) &=& Rw(t) + \sum_{i = 0}^{h-1}  M_i \begin{bmatrix}
				u^{(i)}(t) \\ y^{(i)}(t)
			\end{bmatrix},  \label{obsb}
		\end{IEEEeqnarray}
	\end{subequations}
where $l,h\in\mathbb{N}_0$, $N\in\R^{l\times l}$, $H \in\R^{l\times (k+p)}$, $R\in \R^{r\times l}$, $M_{i}\in\R^{r\times (k+p)}$, $i=0,\ldots,h-1$. A system of the form~\eqref{obs} will be called generalized functional observer because it includes derivatives of the input and output variables. The integers~$l$ and~$h$ denote the order and index of the observer, respectively.
	
We exploit the behavior $\mathscr{B}$ to give a precise definition of generalized functional observers for~\eqref{dls}, similar to~\cite[Def.~3.2]{berger2019disturbance}.
	\begin{definition}\label{def:observer}
		System \eqref{obs} is said to be a generalized functional observer for \eqref{dls}, if for every  $(x,u,y,z) \in \mathscr{B}$ there exists $w\in \mathcal{AC}_{\loc}(\R;\R^l)$ and $\hat z\in \mathcal{L}^1_{\loc}(\R;\R^r)$ such that $(w,u,y,\hat z)$ satisfy~\eqref{obs} for almost all $t\in\R$, and for all $w,\hat z$ with this property we have
		\begin{enumerate}
			\item[(a)] $ \hat{z}(t) - z(t)\to 0$ for $t\to\infty$,
			\item[(b)] if $\hat{z}(0) = z(0)$, then $\hat{z} \myae z$.
		\end{enumerate}
	\end{definition}

Note that the observer design~\eqref{obs} presents a special subclass of the systems considered in~\cite{berger2019disturbance}, where the notion ``partial state observer'' is used and these are constructed as descriptor systems using a system copy. System~\eqref{obs} essentially constitutes a regular descriptor system where the algebraic constraints have been resolved, and results from fixing the dimension of the innovations in the observer class in~\cite{berger2019disturbance}.
	
Now, we establish a relation between partial detectability and the the existence of generalized functional observers for~\eqref{dls}. First, we show that partial detectability is a necessary condition for the existence, essentially corresponding to condition~(a) in Definition~\ref{def:observer}.

\begin{theorem}\label{thm:partdet-nec-obs}
If there exists a generalized functional observer of the form \eqref{obs} for a given system~\eqref{dls}, then~\eqref{dls} is partially detectable.
\end{theorem}
\begin{proof}
Let $(x,0,0,Kx) \in \mathscr{B}$ be arbitrary. Then $w=0$ and $\hat z=0$ satisfy~\eqref{obs} with $u=0$ and $y=0$, thus we have (since~\eqref{obs} is a generalized functional observer)
\begin{equation}\label{eq:necess1}
      z(t) -\hat z(t) = Kx(t) \to 0 \text{ for } t\to\infty.
\end{equation}
Therefore, since in particular $\mathcal{E} \dot x(t) =  \mathcal{A} x(t)$, it follows from~\eqref{eq:necess1} and~\cite[Lem.~A.4]{berger2019disturbance} that condition~(c) in Theorem \ref{thm:kdetectability} is satisfied and hence~\eqref{dls} is partially detectable.
\end{proof}

We like to note that a necessary and sufficient condition for the existence of generalized functional observers in a more general form (which are descriptor systems again) was derived in~\cite[Thm.~3.5]{berger2019disturbance}. Since we consider a smaller class of admissible observers here, the conditions from~\cite{berger2019disturbance} are not sufficient for the existence of an observer anymore. Likewise, although partial detectability is necessary for existence, it is not sufficient. In the remainder of this section, we derive a condition which together with partial detectability also yields the existence of a generalized functional observer. To this end, we also provide an explicit observer design algorithm as follows.

\textbf{Observer design procedure:} Assume that the system~\eqref{dls} is partially detectable.	
First, utilizing Lemma~\ref{pre:lm1}, we compute nonsingular matrices $P_1$ and $Q_1$ such that the matrix pencil $P_1(\lambda \mathcal{E} - \mathcal{A})Q_1$ is in QKF \eqref{eq:qkf}. Then, we compute a non-singular matrix $U_1$ such that \begin{equation}\label{eq:U1}
    U_1^{-1} J_fU_1 = \begin{bmatrix}
	J_{f_1} &  \\ & J_{f_2}
	\end{bmatrix},
\end{equation} where $\sigma(J_{f_1}) \subseteq \bar{\mathbb{C}}^+$ and $\sigma(J_{f_2}) \subseteq \mathbb{C}^-$. The existence of such $U_1$ is guaranteed by the Jordan canonical form of $J_f$. Further, we find (e.g. using the SVD or QR factorization) a non-singular matrix $U_2$ such that
	\begin{equation*}
	U_2 (\lambda E_{\eta} - A_{\eta}) = \begin{bmatrix} \lambda I_{n_\eta} - A_{\eta,1}  \\ -A_{\eta,2} \end{bmatrix}.
	\end{equation*}
Finally, we define
	\begin{IEEEeqnarray}{rCl}
		&P := \blkdiag(I,U_1^{-1},I,U_2)P_1, ~Q :=Q_1\blkdiag(I,U_1,I,I) ,& \nonumber \\ & P\begin{bmatrix} B & 0 \\ D & -I_p \end{bmatrix} = \begin{bmatrix}
			B_{\epsilon} \\ B_{f_1}\\ B_{f_2}\\ B_{\sigma} \\ B_{\eta_1}\\ B_{\eta_2}
		\end{bmatrix},\ KQ = \begin{bmatrix}
			K_{\epsilon}^\top  \\ K_{f_1}^\top \\ K_{f_2}^\top \\ K_{\sigma}^\top \\ K_{\eta}^\top
		\end{bmatrix}^\top, \text{ and }\, x = Q\begin{bmatrix}
			x_{\epsilon}\\ x_{f_1}\\ x_{f_2}\\ x_{\sigma} \\ x_{\eta}
		\end{bmatrix}\!. &\nonumber
		\label{eq:KQ}
	\end{IEEEeqnarray}
Since system~\eqref{dls} is partially detectable, it follows from Lemma~\ref{lm:obsv} that $K_{\epsilon} = 0$ and $K_{f_1} = 0$. Thus, in the new coordinates, the problem of generalized functional observer design for system~\eqref{dls} reduces to the problem of generalized functional observer design for
\begin{subequations}\label{dls5}
\begin{IEEEeqnarray}{rCl}
J_{\sigma} \dot{x}_{\sigma}(t) &=& x_{\sigma(t)} + B_{\sigma}\bar{u}(t), \label{dls5a} \\
\begin{bmatrix}
\dot{x}_{f_2}(t) \\ \dot{x}_{\eta}(t) \end{bmatrix} &=& \begin{bmatrix}
		J_{f_2} \\ & A_{\eta_1}
			\end{bmatrix} \begin{bmatrix}
				x_{f_2}(t) \\ x_{\eta}(t)
			\end{bmatrix} + \begin{bmatrix}
				B_{f_2} \\ B_{\eta_1}
			\end{bmatrix}\bar{u}(t), \label{dls5b}\\
			0  &=& A_{\eta_2} x_{\eta}(t) + B_{\eta_2}\bar{u}(t) , \label{dls5d} \\
			z(t) &=&  K_{\sigma}x_{\sigma}(t) + \begin{bmatrix}
				K_{f_2} & K_{\eta}
			\end{bmatrix} \begin{bmatrix}
				x_{f_2}(t) \\ x_{\eta}(t)
			\end{bmatrix}, \label{dls5e}
		\end{IEEEeqnarray}
	\end{subequations}
	where $\bar u= \begin{bmatrix}
				u \\ y
	\end{bmatrix}$. Since $\rank \begin{bmatrix} \lambda I_{n_\eta} - A_{\eta,1}  \\ -A_{\eta,2} \end{bmatrix}=n_\eta$ for all $\lambda \in\mathbb{C}$ by Lemma \ref{pre:lm1}, there exists $L\in\R^{n_\eta \times (m_\eta - n_\eta)}$ such that $\sigma(A_{\eta,1} - L A_{\eta,2}) \subseteq \mathbb{C}^-$. Define $N := \begin{bmatrix} J_{f_2} & 0 \\ 0 & A_{\eta,1} - L A_{\eta,2} \end{bmatrix}$ and $R := \begin{bmatrix} K_{f_2} & K_{\eta}	\end{bmatrix}$.

In the next theorem, we derive an additional condition on $N$ and $R$ under which the following system is a generalized functional observer for~\eqref{dls}:
\begin{subequations}\label{obs2}
	\begin{IEEEeqnarray}{rCl}
	\dot{w}(t) &=&  Nw(t) + \begin{bmatrix} B_{f_2}  \\ B_{\eta,1} - L B_{\eta,2} \end{bmatrix} \bar{u}(t),\qquad\ \label{obs2a} \\
	\hat{z}(t) &=& Rw(t) - \sum_{i = 0}^{h-1} K_{\sigma} J_{\sigma}^i B_{\sigma}\bar{u}^{(i)}(t) , \label{obs2b}
		\end{IEEEeqnarray}
	\end{subequations}
where $h$ is the index of nilpotency of $J_\sigma$. Before stating the result, we comment on the observer design procedure.

\begin{remark}\label{rem:stable:alg}
The above described observer design procedure can actually be turned into a numerically stable algorithm by
\begin{itemize}
  \item using the staircase form from \cite{beelen1988improved} instead of the QKF,
  \item using the GUPTRI form from \cite{demmel1993generalized} for the decomposition of the stable and unstable eigenvalues of the matrix $J_f$ in \eqref{eq:U1},
\end{itemize}
and some other small but straightforward modifications. However, the presentation of this algorithm is quite technical and we leave the details to the reader.
\end{remark}

\begin{theorem}\label{thm:obsv}
Consider a system \eqref{dls} which is partially detectable and consider  an observer candidate of the form \eqref{obs2} with gain matrix $L\in\R^{n_\eta \times (m_\eta - n_\eta)}$ such that $\sigma(A_{\eta,1} - L A_{\eta,2}) \subseteq \mathbb{C}^-$. If
\begin{equation}\label{thm:obsv:eq1}
    \rank R = \rank \mathcal{O}(R,N),
\end{equation}
where $l = n_{f_2} + n_{\eta}$ and \\ $\mathcal{O}(R,N) := \begin{bmatrix}
R^\top & (RN)^\top & (RN^2)^\top & \ldots & (RN^{l-1})^\top \end{bmatrix}^\top$,
then \eqref{obs2} is a generalized functional observer for \eqref{dls}.
\end{theorem}

\begin{proof}
By construction we have that $\sigma(N) \subseteq \mathbb{C}^-$. Converting system~\eqref{dls} into the form \eqref{dls5} by the transformation described in the observer design procedure, we find that for any $(x,u,y,z)\in\mathcal{B}$ we have that $(x_\sigma, x_{f_2}, x_\eta, \bar u, z)$ is a solution of~\eqref{dls5}. Clearly, $\bar u$ must be sufficiently smooth by~\eqref{dls5a} and hence for all $(x,u,y,z)\in\mathcal{B}$ there exists a solution $(w, \hat z)$ of \eqref{obs2}. For any such solution, we define $e(t) =  z(t) - \hat z(t)$ and $e_1(t) = \begin{bmatrix}
		x_{f_2}(t) \\ x_{\eta}(t)
			\end{bmatrix} - w(t)$,
then using~\eqref{dls5d} we obtain
\begin{subequations}\label{error}
    \begin{IEEEeqnarray}{rCl}
    \dot e_1(t) &=& Ne_1(t) \!+\!\! \begin{bmatrix} 0 \\  L \big( \underset{=0}{\underbrace{A_{\eta,2} x_{\eta}(t) \!+\!   B_{\eta_2}\bar{u}(t)}}\big) \end{bmatrix}\! = Ne_1(t),\quad\label{error1} \\
    e(t) &=& Re_1(t). \label{error2}
\end{IEEEeqnarray}
\end{subequations}
Since $\sigma(N) \subseteq \mathbb{C}^-$, $e(t)\to 0$ for $t\to\infty$. Therefore $\hat{z}(t) \to z(t)$ for $t \to \infty$. Furthermore, if $z(0) = \hat z(0)$, then $e_1(0) \in \ker R $. Since $e_1$ satisfies~\eqref{error1}, it is a direct consequence of condition~\eqref{thm:obsv:eq1} and~\cite[Lem.~A.1]{berger2019disturbance} that $e(t) = Re_1(t) = 0$ for all $t \ge 0$, \emph{i.e.},  $\hat{z} \myae z$.
\end{proof}

\begin{remark}\label{rem:obs}
Denote by $\Sigma_1$ the class of descriptor systems~\eqref{dls} which are partially detectable, by $\Sigma_2$ the class for which there exists a generalized functional observer of the form~\eqref{obs}, and by $\Sigma_3$ the class for which there exists a gain matrix $L$ with $\sigma(N)\subseteq \mathbb{C}^-$ and which satisfies~\eqref{thm:obsv:eq1}. By Theorem~\ref{thm:obsv} the set $\Sigma_3$ consists of those descriptor systems for which~\eqref{obs2} is a generalized functional observer. Furthermore, Theorems~\ref{thm:partdet-nec-obs} and~\ref{thm:obsv} imply the following inclusion:
\[
    \Sigma_1 \supseteq \Sigma_2 \supseteq \Sigma_3.
\]
Obviously, $\Sigma_1 \neq \Sigma_2$ and to see that also $\Sigma_2 \neq \Sigma_3$ consider the following example for a system~\eqref{dls}:
\[
    E=\begin{bmatrix} 1&0\\ 0&1\end{bmatrix},\quad A = \begin{bmatrix} 0&0\\ 1&0\end{bmatrix},\quad B =  \begin{bmatrix} 1\\ 1\end{bmatrix},
\]
\[
    C = \begin{bmatrix} 0 & 1\end{bmatrix},\quad D = 1,\quad K = \begin{bmatrix} 1 & 0\end{bmatrix}.
\]
Then with $P=I_3$ and $Q=I_2$ the system is already in the form~\eqref{dls5} given by
\begin{align*}
    \dot x_\eta(t) &= \begin{bmatrix} 0&0\\ 1&0\end{bmatrix} x_\eta(t) \begin{bmatrix} 1&0\\ 1&0\end{bmatrix} \bar u(t),\\
    0 &= \begin{bmatrix} 0 & 1\end{bmatrix} x_\eta(t) + \begin{bmatrix} 1 & -1\end{bmatrix} \bar u(t),\\
    z(t) &= \begin{bmatrix} 1 & 0\end{bmatrix} x_\eta(t).
\end{align*}
Resolving the equations we find that
\[
    z(t) = \dot y(t) - \dot u(t) - u(t),
\]
and hence there clearly exists a generalized functional observer of the form~\eqref{obs} for the system. However, any observer of the form~\eqref{obs2} with $L=\begin{smallbmatrix} \alpha\\ \beta\end{smallbmatrix}$ reads
\begin{align*}
    \dot w(t) &= \begin{bmatrix} 0&-\alpha\\ 1&-\beta\end{bmatrix} w(t) + \begin{bmatrix} 1-\alpha &\alpha\\ 1-\beta&\beta\end{bmatrix} \bar u(t),\\
    \hat z(t) &= \begin{bmatrix} 1 & 0\end{bmatrix} w(t)
\end{align*}
and clearly $\hat z$ is independent of the derivatives of $y$ and $u$, so it cannot be an estimate for~$z$. And indeed, to have $\sigma(N)\subseteq \mathbb{C}^-$ requires $\alpha\neq 0$, but then
\[
    \rank \begin{bmatrix} 1 & 0\end{bmatrix} = 1 \neq 2 = \rank \begin{bmatrix}1&0\\ 0&-\alpha\end{bmatrix},
\]
so~\eqref{thm:obsv:eq1} is not satisfied. Thus, the system belongs to $\Sigma_2$, but not to $\Sigma_3$.
\end{remark}

\begin{remark}\label{rem4:estimator}
In a couple of works, e.g. \cite{jaiswal2021functional,darouach2022functional,darouach2012functional,darouach2017functional}, an observer is only defined by property (a) in Definition \ref{def:observer} and property (b) is not required. In this case the observer is essentially only an asymptotic estimator. Notably, if we give up condition (b) in Definition \ref{def:observer}, then partial detectability of system \eqref{dls} is indeed a necessary and sufficient condition for the existence of a generalized functional observer, which can be easily deduced from the proofs of Theorems \ref{thm:partdet-nec-obs} and \ref{thm:obsv}.
\end{remark}	
	
\section{Numerical illustration}\label{numerical}
In this section, we illustrate our theoretical findings and validate the proposed algorithm with a numerical example. Consider system \eqref{dls} with the coefficient matrices
	\begin{IEEEeqnarray*}{rCl}
		&E = \begin{bmatrix}
			0 & 0 & 0 & 0 & 0 \\
			0 & 0 & 0 & -1 & 0 \\
			1 & 0 & 0 & 0 & 0 \\
			0 & 0 & -1 & 0 & 0 \\
		    0 & 0 & 0 & 0 & 1 \\
	        0 & 0 & 0 & 0 & 0
    \end{bmatrix},~
		A = \begin{bmatrix}
			0 & 0 & -1 & 0 & 0 \\
			2 & 0 & 0 & 1 & 0 \\
			1 & 0 & 0 & 0 & 0 \\
			0 & -1 & 0 & 0 & 0 \\
			0 & 0 & 0 & 0 & -1 \\
			0 & 0 & 0 & 0 & 1
		\end{bmatrix},& \\
	& B = \begin{bmatrix} -1 \\ 0 \\ 0 \\ 0 \\ 1 \\ 0 \end{bmatrix},\ 	C = \begin{bmatrix}
			0 \\ 0 \\ 1 \\ 0 \\ 0 \end{bmatrix}^\top, ~D = 0, \text{ and }
		K = \begin{bmatrix}
			-1 \\ 1 \\ 1 \\ -1 \\ 1 \end{bmatrix}^\top.  &
	\end{IEEEeqnarray*}
This system is partially detectable by Theorem \ref{thm:kdetectability}. Following the observer design procedure in Section \ref{obsvdesign}, we obtain
	\begin{IEEEeqnarray*}{rCl}
	 P &=& \begin{bmatrix}
	0 & 0 & -1 & 0 & 0 & 0 & 0 \\
	0 & -1 & 1 & 0 & 0 & 0 & 0 \\
	0 & 0 & 0 & -1 & 0 & 0 & 0 \\
	-0.5 & 0 & 0 & 0 & 0 & 0 & 0.5 \\
	0 & 0 & 0 & 0 & 1 & 0 & 0 \\
	0 & 0 & 0 & 0 & 0 & 1 & 0 \\
	 0.5 & 0 & 0 & 0 & 0 & 0 & 0.5 \end{bmatrix}  \\
  \text{ and }
	Q &=& \begin{bmatrix}
	-1 & 0 & 0 & 0 & 0 \\
	0 & 0 & 1 & 0 & 0 \\
	0 & 0 & 0 & 1 & 0 \\
	1 & 1 & 0 & 0 & 0 \\
	0 & 0 & 0 & 0 & 1
	\end{bmatrix}.
\end{IEEEeqnarray*}
Furthermore, the reduced system \eqref{dls5} is of the form
	\begin{IEEEeqnarray*}{rCl}
		\begin{bmatrix} 0 & 1 \\ 0 & 0 \end{bmatrix}\dot{x}_{\sigma}(t) &=& x_{\sigma}(t) + \begin{bmatrix}  0 & 0 \\ 0.5 & -0.5  \end{bmatrix}\bar{u}(t) \\
		\begin{bmatrix}
			\dot{x}_{f_2}(t) \\ \dot{x}_{\eta}(t)
		\end{bmatrix} &=& \begin{bmatrix}  -1 & 0 \\ 0 & -1 \end{bmatrix}
	\begin{bmatrix}
		x_{f_2}(t) \\ x_{\eta}(t)
	\end{bmatrix} + \begin{bmatrix}
	 0 & 0 \\ 1 & 0 \end{bmatrix}\bar{u}(t) \\
		0 &=& \begin{bmatrix} 1 \\ 0 \end{bmatrix}x_{\eta}(t) + \begin{bmatrix}
0 & 0 \\ 0.5 & 0.5
\end{bmatrix}\bar{u}(t) \\
		z(t) &=& \begin{bmatrix}
			  -1 &  1 \end{bmatrix}\begin{bmatrix}
			 x_{f_2}(t) \\ x_{\eta}(t)
		 \end{bmatrix} + \begin{bmatrix}
			1 & 1 \end{bmatrix}x_{\sigma}(t).
	\end{IEEEeqnarray*}
Since the matrix $\begin{bmatrix}
    J_{f_2} & 0\\ 0 & A_{\eta,1}
\end{bmatrix}$ is stable, by choosing $L = 0$ we find that condition \eqref{thm:obsv:eq1} is satisfied and we obtain a generalized functional observer of the form~\eqref{obs} given by
\begin{IEEEeqnarray*}{rCl}
			\dot{w}(t) &=& \begin{bmatrix}
			-1 & 0 \\ 0 & -1
		\end{bmatrix}w(t) + \begin{bmatrix}
		0 & 0 \\ 1 & 0
	\end{bmatrix} \bar{u}(t), \\
		\hat{z}(t) &=& \begin{bmatrix}
			-1 & 1
		\end{bmatrix}w(t) + \begin{bmatrix} 0.5 & -0.5 \end{bmatrix}\bar{u}(t) + \begin{bmatrix} 0.5 & -0.5 \end{bmatrix} \bar{u}^{(1)}(t).
	\end{IEEEeqnarray*}

	\begin{figure}
		\centering
		\begin{center}
			\includegraphics[scale=0.5]{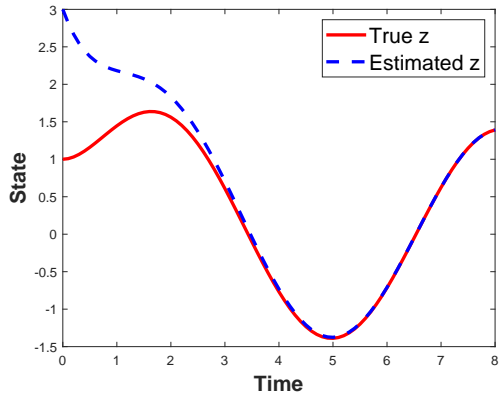}
			\caption{Time response of true and estimated $z(t)$.}
			\label{fig:fig1}
		\end{center}
	\end{figure}
	\begin{figure}
		\centering
		\begin{center}
			\includegraphics[scale=0.5]{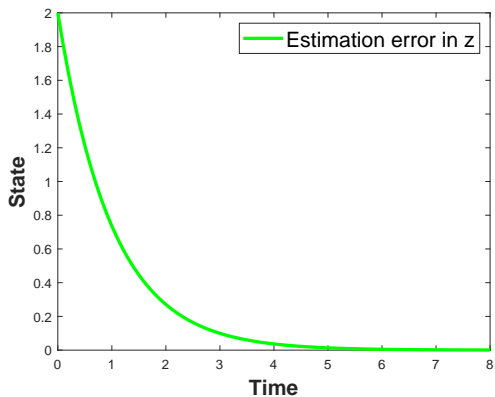}
			\caption{Time response of observer error $e(t)$.}
			\label{fig:fig2}
		\end{center}
	\end{figure}
	Figure \ref{fig:fig1} and \ref{fig:fig2} show the responses of the true and estimated values
	of~$z(t)$ and the time response of the estimation error, respectively. The simulation is realized in MATLAB R$2020$a and the reduced system~\eqref{dls5} and the observer~\eqref{obs2} are solved
by using the ode$15$s solver with relative tolerance $10^{-6}$.	For the simulation, we used the parameters $\begin{bmatrix}
	    x_{f_2}(0) \\ x_{\eta}(0)
	\end{bmatrix} = \begin{bmatrix}
	    1\\2
	\end{bmatrix}$, $w(0) =\begin{bmatrix}
	    3 \\6
	\end{bmatrix}$, and $\bar{u}(t) = \begin{bmatrix}
	    \sin(t) \\ \cos(t)
	\end{bmatrix}$.

\section{Conclusion}\label{conc}
In the present paper, we established a precise mathematical definition and the algebraic characterization of partial detectability for LTI descriptor systems. In addition to this, incorrect results for the algebraic characterization of partial detectability of state space systems in the literature have been pointed out and corrected. It turned out that the proposed concept of partial detectability is a necessary  condition for the existence of generalized functional observers for LTI descriptor systems. Together with a new rank condition, it is also sufficient and a step-by-step observer design procedure was presented. Future research will focus on closing the gap between the system classes discussed in Remark \ref{rem:obs} by additional conditions and finding a characterization of descriptor systems for which a generalized functional observer exists.

\bibliography{funbibfile}	
	
\end{document}